\documentclass[12pt]{elsarticle}
\usepackage[cp1251]{inputenc}
\usepackage[english]{babel}
\usepackage{amsmath}
\usepackage{amssymb}
\usepackage{amsfonts}
\usepackage{mathtools}
\usepackage{array}
\usepackage{rotating}
\usepackage{graphicx}
\usepackage{floatflt,epsfig}
\usepackage{lineno,hyperref}
\usepackage{enumerate}
\usepackage{colortbl}
\usepackage{array,tabularx,tabulary,booktabs}
\usepackage{longtable}
\usepackage{multirow}
\usepackage{wrapfig}
\usepackage{subcaption}
\usepackage[table]{xcolor}
\usepackage{listings}
\usepackage{url}
\usepackage{hyperref}

\journal{The Art of Discrete and Applied Mathematics}
\newtheorem{lemma}{Lemma}

\newtheorem{theorem}{Theorem}

\newcommand{\proof}{\medskip\noindent{\bf Proof.~}}
\begin{document}
\renewcommand{\abstractname}{Abstract}
\renewcommand{\refname}{References}
\renewcommand{\tablename}{Figure.}
\renewcommand{\arraystretch}{0.9}
\thispagestyle{empty}
\sloppy

\begin{frontmatter}
\title{On strictly Deza graphs derived from the Berlekamp-Van Lint-Seidel graph}

\author[01]{Soe~Soe~Zaw}
\ead{soesoezaw313@gmail.com}

\address[01]{Shanghai Jiao Tong University, 800 Dongchuan RD., 200240, 
Minhang District, Shanghai, China}

\begin{abstract}
In this paper, we find strictly Deza graphs that can be obtained from the Berlekamp-Van Lint-Seidel graph by dual Seidel switching. 

\end{abstract}

\begin{keyword}
dual Seidel switching, Berlekamp-Van Lint-Seidel graph, divisible design graph, Deza graph
\vspace{\baselineskip}
\MSC[2010] 05C25\sep 05E30
\end{keyword}
\end{frontmatter}

\section{Introduction}\label{sec0}

Goryainov et al. \cite{GHKS19} gave a characterisation of strictly Deza graphs with parameters $(n,k,k-1,a)$ and $\beta = 1$. They found that such strictly Deza graphs necessarily come from strongly regular graphs having the property $\lambda - \mu = -1$ and can be obtained via two operations: strong product with an edge and the dual Seidel switching \cite{H84}. We are still far away from getting a  classification of strongly regular graphs with $\lambda - \mu = -1$ \cite{AJMP94}. 

It is known that if $\lambda = 0$ and $\mu = 1$, then such a strongly regular graph is either the pentagon, or the Petersen graph, or the Hoffman-Singleton graph, or a hypothetical strongly regular graph with parameters $(3250,57,0,1)$. 

Berlekamp et al. studied strongly regular graphs with $\lambda = 1$ and $\mu = 2$ \cite{BLS73}. It was shown that such a strongly regular graph has parameters either $(9,4,1,2)$ (the only such a graph is $3\times 3$-lattice), or $(99,14,1,2)$, or $(243,22,1,2)$, or $(6273,112,1,2)$, or $(494019,994,1,2)$. Berlekamp et al. further constructed a graph with parameters $(243,22,1,2)$, which is known as the Berlekamp-Van Lint-Seidel graph, but its uniqueness as well as the existence of graphs for the other three parameter tuples remain undecided. In particular, for the tuple $(99,14,1,2)$, this problem is known as the Conway's 99-graph problem.

The smallest feasible parameter tuples of strongly regular graphs with $\lambda = 2$, $\mu = 3$ and $\lambda = 3$, $\mu = 4$ are $(364,33,2,3)$, $(676,45,2,3)$ and $(901,60,3,4)$, respectively   \cite{BDB}, and it is unknown if strongly regular graphs with such parameters exist.

In \cite{GHKS19}, some examples of strictly Deza graphs with parameters $(n,k,k-1,a)$ and $\beta = 1$
were given. In particular, dual Seidel switching was applied to the Petersen graph, the Hoffman-Singleton graph, Paley graphs of square order. In this paper, we investigate if dual Seidel switching can be applied to the Berlekamp-Van Lint-Seidel graph or its complement.

\section{Preliminaries}
We consider undirected graphs without loops or multiple edges.

A $k$-regular graph $\Gamma$ on $v$ vertices is called \emph{strongly regular} with parameters $(v,k,\lambda,\mu)$, $0<k<v-1$,
if any two distinct vertices $x,y$ in $\Gamma$ have $\lambda$ common neighbours when $x,y$ are adjacent and $\mu$ common neighbours if $x,y$ are non-adjacent. For a vertex $x$ in a graph $\Gamma$,
the \emph{neighbourhood} $\Gamma(x)$ is the set of all neighbours of $x$ in $\Gamma$.

\begin{lemma}[\cite{BCN89}, Theorem 1.3.1(i)]\label{SRGSpectrum}
Let $\Gamma$ be a strongly regular graph with parameters $(v,k,\lambda,\mu)$, $\mu \ne 0$, $\mu \ne k$.
Then $\Gamma$ has three distinct eigenvalues $k, r, s$, where $k > r > 0 > s$ and the eigenvalues
$r,s$ satisfy the quadratic equation $x^2 + (\mu - \lambda)x + (\mu-k) = 0$.
\end{lemma}

For a graph $\Gamma$, denote by $\overline{\Gamma}$ the complement of $\Gamma$.
\begin{lemma}[\cite{BCN89}, Theorem 1.3.1(vi)]\label{SRGComplement}
Let $\Gamma$ be a strongly regular graph with parameters $(v,k,\lambda,\mu)$. Then the complement $\overline{\Gamma}$ of
$\Gamma$ is a strongly regular graph with parameters $(v,v-k-1, v-2k+\mu-2, v-2k+\lambda)$
and eigenvalues $v-k-1, -s-1,-r-1$.
\end{lemma}

A $k$-regular graph $\Delta$ on $v$ vertices is called a \emph{Deza graph} with parameters $(v,k,b,a), b \ge a$, if the number of common neighbours of any two distinct vertices in $\Delta$ takes on the two values $a$ or $b$. A Deza graph $\Delta$ is called a \emph{strictly Deza graph},
if it has diameter $2$ and is not strongly regular.
The following lemma gives a construction of strictly Deza graphs, which is known as \emph{dual Seidel switching}.

\begin{lemma}[\cite{EFHHH99}, Theorem 3.1]\label{Construction}
Let $\Gamma$ be a strongly regular graph with parameters $(v,k,\lambda,\mu)$, $k\ne\mu$, $\lambda\ne\mu$ and adjacency matrix $M$. Let $P$ be a permutation matrix that represents an involution $\phi$ of $\Gamma$ that interchanges only non-adjacent vertices. Then $PM$ is the adjacency matrix of a strictly Deza graph $\Delta$
with parameters $(v,k,b,a)$, where $b = max(\lambda,\mu)$ and $a = min(\lambda,\mu)$.
\end{lemma}
Since $\phi$ in Lemma \ref{Construction} represents an involution,
the matrix $PM$ is obtained from the matrix $M$ by a permutation of rows in all pairs of rows with indexes $i_1$ and $i_2$, such that $\phi(i_1) = i_2$ and $\phi(i_2) = i_1$. Lemma \ref{Images} follows immediately from Lemma \ref{Construction} and shows what is the neighbourhood of a vertex of the graph $\Delta$.
\begin{lemma}\label{Images} For the neighbourhood $\Delta(u)$ of a vertex $u$ of the graph $\Delta$ from Lemma \ref{Construction}, the following conditions hold:
$$
\Delta(u) =
\left\{
  \begin{array}{ll}
    \Gamma(u), & \hbox{if $\phi(u) = u$;} \\
    \Gamma(\phi(u)), & \hbox{if $\phi(u) \ne u$.}
  \end{array}
\right.
$$
\end{lemma}

In \cite[Theorem 2]{GHKS19}, it was shown that the strong product with an edge and dual Seidel switching is the only method to obtain strictly Deza graphs with $k = b + 1$.
Recall that the graph \emph{strong product} of two graphs $\Gamma_1$ and $\Gamma_2$ has vertex set $V(G_1)\times V(G_2)$ and two distinct vertices  $(v_1,v_2)$ and $(u_1,u_2)$ are connected iff they are adjacent or equal in each coordinate, i.e., for $i \in {1,2}$, either $v_i=u_i$ or $\{v_i,u_i\}$ in $E(\Gamma_i)$, where $E(\Gamma_i)$ is the edge set of $\Gamma_i$   \cite{BW04}.

It follows from Lemma \ref{SRGComplement} that, if a strongly regular graph $\Gamma$ has the property $\lambda - \mu = -1$, then the complementary graph $\overline{\Gamma}$ has the property $\overline{\lambda} - \overline{\mu} = -1$ as well. Thus, according to \cite[Theorem 2]{GHKS19}, we concentrate on order 2 automorphisms of $\Gamma$ that interchange either only non-adjacent vertices or only adjacent vertices. 

Let $G$ be a group and $S$ be an inverse-closed identity-free subset in $G$. The graph on $G$ with two vertices $x,y$ being adjacent whenever $xy^{-1}$ belongs to $S$ is called the \emph{Cayley graph} of the group $G$ with \emph{connection set} $S$ and is denoted by $Cay(G,S)$.

\section{The Berlekamp-Van Lint-Seidel graph and dual Seidel switching}
The \emph{Berlekamp-Van Lint-Seidel graph}, denoted by $\Gamma$, is the coset graph of the ternary Golay code \cite[Section 11.3]{BCN89}. This graph is known to be strongly regular with parameters $(243,22,1,2)$.

In this section, we deal with two more ways to define this graph and give a description of the involutions of $\Gamma$ and $\overline{\Gamma}$ suitable for dual Seidel switching.

The main result of this paper is the following theorem.

\begin{theorem}\label{main} The following statements hold.\\
{\rm (1)} $\Gamma$ has no order 2 automorphisms that interchange only adjacent vertices;\\
{\rm (2)} $\Gamma$ has the unique (up to conjugation) order 2 automorphism
that interchanges only non-adjacent vertices.
\end{theorem}

To prove Theorem \ref{main}, we prove two lemmas, which imply the truth of the theorem statements.

\subsection{$\Gamma$ from the Mathieu group $M_{11}$}
By ATLAS of Group Representations the
Mathieu group $M_{11}$ can be represented \cite{W} by $5\times5$ matrices over $GF(3)$ as follows. 
Put
$$
x := 
\begin{pmatrix}
 0 & 2 & 1 & 0 & 0 \\
 2 & 1 & 1 & 2 & 2 \\
 0 & 1 & 1 & 2 & 2 \\
 1 & 0 & 2 & 2 & 1 \\
 1 & 2 & 2 & 2 & 0
\end{pmatrix},~
y :=
\begin{pmatrix}
 0 & 0 & 2 & 0 & 2 \\
 1 & 1 & 2 & 2 & 0 \\ 
 2 & 2 & 2 & 2 & 2 \\
 1 & 2 & 1 & 1 & 0 \\
 2 & 2 & 0 & 2 & 1
\end{pmatrix}.
$$
Then the group $G:=\langle x,y \rangle$ is isomorphic to $M_{11}$, where $x$ is an involution. Let $V(5,3)$ denote the 5-dimensional vector space of over $GF(3)$. We regard the elements of $V(5,3)$ as rows and consider the action of $G$ on $V(5,3)$ by the right multiplication, which has two orbits of size $22$ and $220$ on the nonzero vectors. The orbit of size $22$ is given by the set $$S_1:=\{
\pm(1,0,0,0,0),
\pm(0,0,1,0,1),
\pm(0,1,0,1,0),
\pm(0,1,2,0,0),$$$$
\pm(0,0,1,2,1),
\pm(0,1,0,1,2),
\pm(1,1,2,0,2),
\pm(1,0,0,1,2),$$$$
\pm(1,0,2,1,0),
\pm(1,1,0,0,2),
\pm(1,1,2,1,0)
\},$$
and, moreover, $\Gamma$ is isomorphic to the Cayley graph $Cay(V(5,3),S_1)$. Since $G$ stabilises $S_1$ setwise, $G$ is a subgroup in the automorphism group of $\Gamma$, which is known (see \cite{B}) to be isomorphic to the group $3^5 : (2 \times M_{11})$. The fact that $M_{11}$ has precisely one class of conjugate involutions implies that the automorphism group of $\Gamma$ has precisely three classes of conjugate involutions. Let $e$ be the identity matrix from $G$. Note that $-e$ does not belong to $G$, but the multiplication by $-e$ is an involution of the automorphism group of $Cay(V(5,3),S_1)$, which means that the three pairwise non-conjugate involutions of the automorphism group of $Cay(V(5,3),S_1)$ are given by the right multiplication by $-e$, $x$ and $-x$.

\begin{lemma}\label{ex}
The following statements hold. \\
(1) The involution $-e$ interchanges adjacent vertices as well as non-adjacent ones;\\
(2) The involution $-x$ interchanges adjacent vertices as well as non-adjacent ones.
\end{lemma}
\proof
(1) This involution fixes the zero vector and moves all non-zero vectors by swapping every two elements that are additive inverses of each other. In the graph  $Cay(V(5,3),S_1)$, two additive inverses are adjacent whenever both of them belong to $S_1$. It means that the involution $-e$ interchanges adjacent vertices as well as non-adjacent ones.

(2) On the one hand, the involution $-x$ swaps the vertices $(0,1,0,1,0)$ and $(0,2,0,2,0)$, which are adjacent in $Cay(V(5,3),S_1)$. On the other hand, $-x$ swaps the vertices $(1,0,0,0,0)$ and $(0,2,1,0,0)$, which are not adjacent in $Cay(V(5,3),S_1)$.
$\square$

\medskip

In view of Lemma \ref{ex}, it remains to check the inner involution $x$. In the next subsection, we explore one more definition of the Berlekamp-Van Lint-Seidel graph and give a very natural description of the involution $x$.

\subsection{Specific parity-check matrix}
Recall that, for a positive integer $n$ and a prime power $q$, $V(n,q)$ denotes the $n$-dimensional vector space over the finite field $\mathbb{F}_q$.
The ternary Golay code can be constructed as the 6-dimensional subspace in $V(11,3)$ consisting of all row-vectors $\textbf{c}$ such that the equality $H\textbf{c}^T = \textbf{0}$ holds, where 
$$
H :=
\begin{bmatrix*}[r]
1&1&1&2&2&0&1&0&0&0&0\\
1&1&2&1&0&2&0&1&0&0&0\\
1&2&1&0&1&2&0&0&1&0&0\\
1&2&0&1&2&1&0&0&0&1&0\\
1&0&2&2&1&1&0&0&0&0&1 
\end{bmatrix*}
$$
is the specific parity check matrix of this code. 
Let $x_{1}, x_{2}, x_{3}, \ldots, x_{11}$ denote the vectors from $V(5, 3)$ that correspond to the columns of $H$. There are 22 vectors of type $\pm x_{i}$ and $220$ vectors of type $ \pm{x_{i}} \pm{x_{j}}$ where $i \neq j;  i,j = 1, 2,\ldots, 11$. The Cayley graph $Cay(V(5,3),S_2)$, where $S_2 := \{\pm x_{1}, \ldots, \pm x_{11}\}$, is known to be isomorphic with the Berlekamp-Van Lint-Seidel graph (see \cite{BLS73}).

\begin{lemma}\label{Compl} The reversal of vectors is an involution of $Cay(V(5,3),S_2)$ that interchanges only non-adjacent vertices.
\end{lemma}
\proof
Obviously, the reversal of vectors is a permutation of the vertex set of $\Gamma$. For a vector $\gamma \in V(5,3)$, denote by $\gamma^r$ the reversed vector. Note that $(S_{2})^r = S_{2}$ holds. Since, for any two vertices $\gamma_1,\gamma_2$ in $\Gamma$, we have $\gamma_1^r - \gamma_2^r = (\gamma_1-\gamma_2)^r$, the reversal is an automorphism of $Cay(V(5,3),S_2)$.

For a vector $(a,b,c,d,e) \in V(5,3)$, consider the difference $(a,b,c,d,e) - (a,b,c,d,e)^r = (a-e,b-d,0,d-b,e-a)$. Note that the first and the fifth coordinates and the second and fourth ones are additive inverses. Since $S_2$ has no such vectors with zero third coordinate, the reversal interchanges only non-adjacent vertices.
$\square$

\section{Concluding remarks}
The following three strictly Deza graphs can be derived from the Berlekamp-Van Lint-Seidel graph $\Gamma$.

First, Lemma \ref{Construction} and Theorem 1(2) give a strictly Deza graph with parameters $(243,22,2,1)$. It has spectrum $\{ 22^1, 5^{48}, 4^{72}, (-4)^{60}, (-5)^{62}\}$ and its automorphism group of order 2592 is a subgroup in the automorhism group of $\Gamma$. 

Further, in view of \cite[Construction 1]{GHKS19}, the strong product $\Gamma[K_{2}]$ of the Berlekamp-Van Lint-Seidel graph with an edge is a strictly Deza graph with parameters $(486,45,44,4)$. It has spectrum
$\{  45^1, 9^{132},  (-1)^{243}, (-9)^{110} \}$.

Finally, the order 2 automorphism from Theorem 1(2) induces an order 2 automorphism of $\Gamma[K_{2}]$ that interchanges only non-adjacent vertices. Applying the dual Seidel switching to $\Gamma[K_{2}]$, we obtain one more strictly Deza graph with parameters $(486,45,44,4)$, which has spectrum $\{45^1, 9^{120}, 1^{108}, (-1)^{135}, (-9)^{122}\}$.

In the connection with the results from \cite{GHKS19}, we point out that both graphs with parameters $(486,45,44,4)$ are divisible design graphs.

\section*{Acknowledgment} \label{Ack}
The author thanks both anonymous referees for their comments and suggestions, which significantly improved the paper. The author thanks Professor Yaokun Wu for his continued support and warm hospitality and Sergey Goryainov for valuable discussions.
This work is supported by NSFC(11671258) and STCSM(17690740800).

\end{document}